\documentclass[10pt,twoside]{article}
\usepackage{Latex-document}
\newtheorem{theo}{Theorem}[section]
\newtheorem{prop}[theo]{Proposition}

\newtheorem{conj}[theo]{Conjecture}

\newcommand{\G}{{\cal G}}

\renewcommand{\thesection}{\arabic{section}}
\renewcommand{\thesubsection}{\thesection.\arabic{subsection}.}
\def\Section#1{\section{\hskip -1em . \hskip 0.6em#1}}
\def\volumeno{I}

\title{\bf Discrete Mathematics: \vskip -2mm Methods and Challenges\vskip 6mm}

\author{Noga Alon\thanks{School of Mathematics and Computer Science,
Raymond and Beverly Sackler Faculty of Exact Sciences, Tel Aviv
University, Tel Aviv, Israel 69978. E-mail: noga@math.tau.ac.il}
\vspace{-0.5cm}}

\date{\vspace{-8mm}}

\begin{document}

\maketitle

\thispagestyle{first} \setcounter{page}{119}

\begin{abstract}

\vskip 3mm

Combinatorics is a fundamental mathematical discipline as well as
an essential component of many mathematical areas,
and its study has
experienced an impressive growth in recent years.
One of the main reasons
for this growth is the tight connection between Discrete Mathematics and
Theoretical Computer Science, and the rapid development of the latter.
While in the past many of the basic combinatorial results were obtained
mainly by ingenuity and detailed reasoning, the modern theory has grown
out of this early stage, and often relies on deep, well developed tools.
This is a survey of two of the main general techniques that played
a crucial role in the development of modern combinatorics; algebraic
methods and probabilistic methods. Both will be illustrated by examples,
focusing on the basic ideas and the connection to other areas.

\vskip 4.5mm


\end{abstract}

\vskip 12mm

\Section{Introduction} \setzero

\vskip-5mm \hspace{5mm}

The originators of the basic concepts of Discrete Mathematics,
the mathematics of finite structures, were
the Hindus, who knew the formulas for the number of permutations of
a set of $n$ elements, and for the number of subsets of cardinality
$k$ in a set of $n$ elements, already in the sixth century. The
beginning of Combinatorics as we know it today started with the work
of Pascal and De Moivre in the $17$th century, and continued in the
$18$th century with the seminal ideas of Euler in Graph Theory, with his
work on partitions and their enumeration, and with his interest in
latin squares. These old results are among the roots of the study
of formal methods of enumeration,
the development of configurations and designs, and the extensive work
on Graph Theory in the last two centuries. The tight connection
between Discrete Mathematics and Theoretical Computer Science,
and the rapid development of the latter in recent years, led to
an increased interest in combinatorial techniques and to
an impressive development of the subject. It also
stimulated the study of algorithmic combinatorics
and combinatorial optimization.

While many of the basic combinatorial results were obtained
mainly by ingenuity and detailed reasoning, without relying
on many deep, well developed tools, the modern theory has
already grown out of this early stage. There are already well
developed enumeration methods, some of which are based on deep
algebraic tools. The probabilistic method initiated by Erd\H{o}s
(and to some extent, by Shannon) became one of the most powerful
tools in the modern theory, and its study has been fruitful to
Combinatorics, as well as to Probability Theory. Algebraic and
topological techniques play a crucial role in the modern theory,
and Polyhedral Combinatorics, Linear Programming and constructions
of designs have been developed extensively. Most of the new
significant results obtained in the area are inevitably based on
the knowledge of these well developed concepts and techniques,
and while there is, of course, still a lot of room for pure ingenuity
in Discrete Mathematics, much of the progress is obtained by relying
on the fast growing accumulated body of knowledge.

Concepts and questions of Discrete Mathematics  appear naturally in
many branches of mathematics, and the area has found applications
in other disciplines as well. These include applications in
Information Theory and Electrical Engineering, in Statistical Physics,
in Chemistry and Molecular Biology, and, of course, in Computer Science.
Combinatorial topics such as Ramsey Theory, Combinatorial Set Theory,
Matroid Theory, Extremal Graph Theory, Combinatorial Geometry
and Discrepancy Theory are related
to a large part of the mathematical and scientific world, and these
topics have already found numerous applications in other fields.
A detailed account of the topics, methods and applications of
Combinatorics can be found in \cite{GGL}.

This paper is mostly a survey of two of the main general
techniques that played a crucial role in the development of modern
combinatorics; algebraic methods and probabilistic methods. Both
will be illustrated by examples, focusing on
the basic ideas and the connection to other areas.
The choice of topics and examples described here is inevitably biased,
and is not meant to be comprehensive.
Yet, it hopefully provides some
of the flavor of the techniques, problems and results in the area in a
way which may be appealing to researchers, even if their main interest
is not Discrete Mathematics.

\Section{Dimension, geometry and information theory} \setzero

\vskip-5mm \hspace{5mm}

Various algebraic techniques have been used successfully in
tackling problems in Discrete Mathematics over the years. These
include tools from Representation Theory  applied extensively
in enumeration problems, spectral techniques used in the study of highly
regular structures, and applications of
properties of polynomials and tools from algebraic geometry
in the theory of Error  Correcting Codes
and in the study of problems in Combinatorial Geometry. These
techniques have numerous interesting applications. Yet, the most
fruitful algebraic technique applied in
combinatorics, which is possibly also the simplest one, is the so-called
dimension argument. In its simplest form, the method can be described
as follows.
In order to bound the cardinality of a discrete
structure $A$,
one maps its elements to vectors in a linear space, and shows that
the set $A$ is mapped to a linearly independent set. It then follows
that the cardinality of $A$ is bounded by the dimension of the
corresponding linear
space. This method is often particularly useful in the solution
of extremal problems in which
the extremal configuration is not unique.  The method is effective
in such cases because bases in a
vector space can be very different from each other and
yet all of them have the same cardinality. Many applications of this
basic idea can be found in \cite{BF}, \cite{Bl}, \cite{GR}.

\subsection{Combinatorial geometry}

\vskip-5mm \hspace{5mm}

An early application of the dimension argument
appears in \cite{LRS}. A set of points $A
\subset R^n$ is a two-distance
set if at most $2$ distinct positive distances are determined by the
points of $A$.  Let $f(n,2)$ denote the
maximum possible size of a two-distance set in
$R^n$.
The set of all $0/1$ vectors in $R^{n+1}$ with exactly two $1$'s
shows that $f(n,2) \geq n(n+1)/2$, and the authors of \cite{LRS}
proved that $f(n,2) \leq (n+1)(n+4)/2$. The upper bound is proved
by associating each point of a two-distance set $A$ with a polynomial
in $n$ variables, and by showing that these polynomials are
linearly independent and all lie in a space of dimension
$(n+1)(n+4)/2$. This has been improved by Blokhuis to
$(n+1)(n+2)/2$, by showing that one can add $n+1$ additional
polynomials that lie in this space
to those obtained from the two-distance set, keeping the
augmented set linearly independent. See \cite{Bl} and its references
for more details. The precise value of $f(n,2)$ is not known.

Borsuk \cite{Bo} asked if any compact set of at least $2$
points in $R^d$ can be
partitioned into at most $d+1$ subsets of smaller diameter.
Let $m(d)$ be the smallest integer $m$ so that any such set
can be partitioned into a most $m$ subsets of smaller diameter.
Borsuk's question is whether $m(d)=d+1$ (the $d+1 $ points
of a simplex show that $m(d)$ is at least $d+1$.)
Kahn and Kalai \cite{KK} gave an example showing that this is not the
case for all sufficiently large $d$,
by applying a theorem of Frankl and Wilson \cite{FW}.
Improved versions of their construction have been obtained by
Nilli in 1994, by Raigorodski  in 1997, by Hinrichs  in 2001
and by Hinrichs and Richter
in 2002.
%%%%%%%%%%%%%
The last two results
are based on some properties of the Leech Lattice and
give a construction showing that
already in dimension $d=298$, more than $d+1$ subsets may be needed.
All the constructions and the proofs of their properties are based on
the dimension argument. Here is a brief sketch of one of them.

Let $n=4p$, where $p$ is an odd prime, and let
${\cal F}$ be the set of all vectors ${\bf x}=(x_1, \ldots ,x_n)
\in \{-1,1\}^n$, where $x_1=1$ and the number of negative coordinates
of ${\bf x}$ is even. One first proves the following.

\begin{equation}
\label{e21}
\mbox{
If ${\cal G} \subset {\cal F}$ contains no two orthogonal vectors
then $|{\cal G}| \leq \sum_{i=0}^{p-1} {n-1 \choose i}.$
}
\end{equation}

This is done by associating each member of $\G$  with a
multilinear polynomial of
degree at most $p-1$ in $n-1$ variables,  so that all the
obtained polynomials are linearly independent.
Having established (\ref{e21}),
define $S=\{{\bf x} *{\bf x}: {\bf x} \in {\cal F} \}$,
where ${\cal F}$ is as above, and
${\bf x}*{\bf x}$
is the tensor product of ${\bf x}$ with itself, i.e., the vector
of length $n^2$, $(x_{ij}:1 \leq i,j \leq n)$, where $x_{ij}=x_i x_j$.
The norm of each vector in $S$
is $n$ and the scalar product between any two members of $S$
is non-negative. Moreover, by
(\ref{e21}) any set of more than $\sum_{i=0}^{p-1}{n-1 \choose i}$
members of $S$ contains an orthogonal pair, i.e., two points
the distance between which is the diameter of $S$. It follows that
$S$ cannot be partitioned into less than
$2^{n-2}/ \sum_{i=0}^{p-1}{n-1 \choose i}$ subsets of smaller diameter.
This shows that $m(d) \geq c_1^{\sqrt d}$ for some $c_1 >1$.
An upper bound of $m(d) \leq c_2^d$ where $c_2 = \sqrt{3/2}+o(1)$
%%%%%%%%%
is known,
but determining the correct order of magnitude of $m(d)$ is an open
question. The following conjecture seems plausible.
\begin{conj}
\label{c91}
There is a constant $c >1$ such that
$ m(d) > c^d $
for all $d \geq 1$.
\end{conj}

An {\em equilateral set} (or a simplex) in a metric space,
is a set $A$, so that the distance between any pair of distinct members
of $A$ is $b$, where $b \neq 0$ is a constant.
Trivially, the maximum cardinality of such
a set in $R^n$ with respect to the (usual) $l_2$-norm
is $n+1$. Somewhat surprisingly, the situation is far more
complicated for the  $l_1$ norms.
The $l_1$-distance between two points $\vec{a}=(a_1, \ldots a_n)$
and $\vec{b}=(b_1, \ldots ,b_n)$ in $R^n$ is
$\|\vec{a}-\vec{b}\|_1=(\sum_{k=1}^n |a_i-b_i|$.
Let $e(l_1^n)$ denote the maximum
possible cardinality of an equilateral set in $l_1^n$.
The set of standard basis vectors and their negatives shows that
$e(l_1^n) \geq 2n$. Kusner \cite{Gu} conjectured that this is
tight, i.e., that $e(l_1^n)=2n$ for all $n$.  For $n \leq 4$
this is proved in \cite{KLS}. For general $n$, the best known
upper bound is $e(l_1^n)\leq c_1 n \log n$ for some absolute
positive constant $c_1$. This is proved in \cite{AP2} by  an
appropriate dimension argument. Each vector in an equilateral
set of $m$ vectors in $R^n$ is mapped to a vector in $l_2^t$
for an appropriate $t=t(m,n)$,
by applying a probabilistic technique involving
randomized rounding. It is then shown, using a simple argument based
on the eigenvalues of the Gram matrix of these new vectors, that they
span a space of dimension at least $c_2 m$, implying that
$c_2 m \leq t(m,n)$ and supplying the desired result. The precise
details require some work, and can be found in \cite{AP2}.

\subsection{Capacities and graph powers}

\vskip-5mm \hspace{5mm}

Let $G=(V,E)$ be a simple, undirected graph.  The power
$G^n$ of $G$  is the graph whose set of vertices is $V^n$
in which
two distinct vertices
$(u_1, u_2, \ldots ,u_n)$ and $(v_1, v_2, \ldots ,v_n)$ are
adjacent iff for all $i$ between $1$ and $n$ either $u_i=v_i$ or
$u_iv_i \in E$.  The {\em Shannon capacity} $c(G)$ of $G$ is the limit
$\lim_{n \rightarrow \infty} (\alpha(G^n))^{1/n}$, where $\alpha(G^n)$
is the maximum size of an independent set of vertices in $G^n$. This
limit exists, by super-multiplicativity, and it is always at least
$\alpha(G)$.

The study of this parameter was introduced by Shannon in \cite{Sh},
motivated by a question in Information Theory. Indeed, if $V$ is the
set of all possible letters a channel can transmit in one use, and two
letters are adjacent if they may be confused, then $\alpha(G^n)$ is
the maximum number of messages that can be transmitted in $n$ uses of
the channel with no danger of confusion.  Thus $c(G)$ represents the
number of distinct messages {\em per use} the channel can communicate
with no error while used many times.

Calculation of $c(G)$ seems to be very hard.  For example $c(C_5) =
\sqrt 5$ was only shown in
1979 by Lov\'asz \cite{Lo}, and $c(C_7)$ remains
unknown.  Certain polynomially computable upper bounds on $c(G)$ are
known including Lov\'asz's theta function $\theta(G)$, and other upper
bounds are due to Haemers and to Schrijver.

Another upper bound, based on the dimension argument and related
to the bound of Haemers \cite{Ha}, is described in \cite{Al21},
where it is applied to solve a problem of Shannon on the capacity
of the disjoint union of two graphs. The {\em (disjoint) union} of
two graphs $G$ and $H$, denoted by $G+H$, is the graph whose
vertex set is the disjoint union of the vertex sets of $G$ and of
$H$ and whose edge set is the (disjoint) union of the edge sets of
$G$ and $H$. If $G$ and $H$ are graphs of two channels, then their
union represents the {\em sum} of the channels corresponding to
the situation where either one of the two channels may be used, a
new choice being made for each transmitted letter. Shannon proved
that for every $G$ and $H$, $c(G+H) \geq c(G) + c(H)$ and that
equality holds in many cases. He conjectured that in fact equality
always holds.  In \cite{Al21} it is shown that this is false in
the following strong sense.

\begin{theo}
\label{t92}
For every $k$ there is a graph $G$ so that the Shannon capacity of
the graph and that of its complement $\overline{G}$ satisfy
$c(G) \leq k, c(\overline{G}) \leq k$, whereas
$c(G+\overline{G}) \geq k^{(1 + o(1))\frac{\log k }{8 \log \log k}}$
and the $o(1)$-term tends to zero as $k$ tends to infinity.
\end{theo}

Therefore, the capacity of the disjoint union of two graphs can be
much bigger than the capacity of each of the two graphs. Strangely
enough, it is not even known if the maximum possible capacity of
a disjoint union of two graphs $G$ and $H$, each of capacity at most
$k$, is bounded  by any function of $k$. It seems very likely that this
is the case.

\Section{Polynomials, addition and graph coloring} \setzero

\vskip-5mm \hspace{5mm}

The study of algebraic varieties, that is, sets of common roots of
systems of polynomials, is the main topic of algebraic geometry.
The most elementary property of a univariate nonzero polynomial
over a field is the fact that it does not have more  roots than
its degree. This elementary property is surprisingly effective in
Combinatorics: it plays a major role in the theory of error
correcting codes, and has many applications in the study of finite
geometries --- see, e.g., \cite{Bl}. A similar property holds for
polynomials of several variables, and can also be used to supply
results in Discrete Mathematics.  In this section we describe a
general result of this type, which is called  in \cite{Al3} {\em
Combinatorial Nullstellensatz}, and briefly sketch some of its
applications in Additive Number Theory and in Graph Theory.

\subsection{Combinatorial nullstellensatz}

\vskip-5mm \hspace{5mm}

Hilbert's  Nullstellensatz (see, e.g., \cite{vdW}) is the
fundamental theorem that asserts that if $F$ is an algebraically
closed field, and $f,g_1, \ldots ,g_m$ are polynomials in
the ring of polynomials $F[x_1, \ldots ,x_n]$, where $f$ vanishes
over all common zeros of $g_1, \ldots ,g_m$, then there is an
integer $k$ and polynomials $h_1, \ldots ,h_m$ in
$F[x_1, \ldots ,x_n]$ so that
$$
f^k=\sum_{i=1}^n h_i g_i.
$$
In the special case $m=n$, where each $g_i$ is a univariate
polynomial of the form $\prod_{s \in S_i} (x_i-s)$ for some
$S_i \subset F$, a
stronger conclusion holds. It can be shown
that if
$F$ is an arbitrary field, $f,g_i, S_i$ are as above,
and $f$ vanishes over all the common zeros of
$g_1, \ldots ,g_n$ (that is; $f(s_1, \ldots ,s_n)=0$
for all $s_i \in S_i$), then there are polynomials
$h_1, \ldots ,h_n \in F[x_1, \ldots ,x_n]$ satisfying
$deg(h_i) \leq deg (f)-deg(g_i)$ so that
$$
f=\sum_{i=1}^n h_i g_i.
$$

As a consequence of the above one can prove the following.
\begin{theo}
\label{t12}
Let $F$ be an arbitrary field, and let $f=f(x_1, \ldots ,x_n)$
be a polynomial in $F[x_1, \ldots ,x_n]$. Suppose the degree
$deg(f)$ of $f$ is $\sum_{i=1}^n t_i$, where each $t_i$ is a
nonnegative integer, and suppose the coefficient of
$\prod_{i=1}^n x_i^{t_i}$ in $f$ is nonzero. If
$S_1, \ldots ,S_n$ are subsets of $F$ with $|S_i|>t_i$,
then there are $s_1 \in S_1, s_2 \in S_2, \ldots, s_n \in S_n$
so that
$$
f(s_1, \ldots ,s_n) \neq 0.
$$
\end{theo}

The detailed proof, as well as many applications,
can be found in \cite{Al3}.
A quick application,
first proved in \cite{AFK}, is the assertion that
for any prime $p$, any loopless graph $G=(V,E)$ with average degree
bigger than $2p-2$ and maximum degree  at most $2p-1$ contains
a $p$-regular subgraph.

To prove it, let $(a_{v,e})_{v \in V, e \in E}$
denote the incidence matrix of $G$
defined by $a_{v,e}=1 $ if $v \in e$ and $a_{v,e}=0$ otherwise.
Associate each edge $e$ of $G$ with a variable $x_e$ and consider
the polynomial
$$
f=\prod_{v \in V} [1-(\sum_{e \in E} a_{v,e} x_e)^{p-1}]
- \prod_{e \in E} (1-x_e),
$$
over $GF(p)$.
Applying Theorem \ref{t12} with $t_i=1$ and $S_i=\{0,1\}$
for all $i$, we conclude that
there are values $x_e \in \{0,1\}$ such that $f(x_e: e \in E) \neq
0$.
It is now easy to check that
in the subgraph consisting of all edges $e \in E$
for which $x_e =1$ all degrees are divisible by $p$, and since
the maximum degree is smaller than $2p$ all positive
degrees are precisely $p$, as needed.

Pyber applied the above result to solve a problem of Erd\H{o}s and
Sauer and prove that any
simple graph on $n$ vertices with at least $200 n \log n$ edges
contains a $3$-regular subgraph.
Pyber, R\"odl and Szemer\'edi  proved that this is not very far from
being best possible, by showing, using
probabilistic arguments,
that there are simple graphs on $n$ vertices with at least
$c n \log \log n$ edges that contain no $3$-regular subgraphs.
See \cite{PRS} for some further related results.

\pagebreak

\subsection{Additive number theory}

\vskip-5mm \hspace{5mm}

The Cauchy-Davenport Theorem, which has numerous
applications in Additive Number Theory, is the statement that
if $p$ is a prime, and $A,B$ are two nonempty subsets of $Z_p$,
then $$|A+B| \geq min\{p,|A|+|B|-1\}.$$

Cauchy proved this theorem in 1813, and applied it to give a
new proof  to a lemma of Lagrange in his well known
1770 paper that shows
that every positive integer is a sum of four squares. Davenport formulated the
theorem as a discrete analogue of a conjecture of Khintchine
about the Schnirelman density of the sum of two sequences of
integers. There are numerous extensions of this result,
see, e.g., \cite{Na}.
A simple algebraic proof of this result is given in
\cite{ANR2},
and its main advantage is that it extends easily and gives
several related results.  This proof can be
described as a simple application of Theorem \ref{t12}.
If $|A|+|B|>p$, then the result is trivial, as the sets $A$ and $g-B$
intersect, for each $g \in Z_p$. Otherwise,
assuming the result is false and $|A + B| \leq |A| +|B|-2$, let
$C$ be a subset
of $Z_p$ satisfying $A+B \subset C$ and $|C|=|A|+|B|-2$. Define
$f=f(x,y)=\prod_{c \in C} (x+y-c)$
and
apply Theorem \ref{t12} with
$t_1=|A|-1,t_2=|B|-1$, $S_1=A, S_2=B$ to get a contradiction.

Using similar (though somewhat more complicated)
arguments, the following related result is proved
in \cite{ANR2}.

\begin{prop}
\label{p93}
Let $p$ be a prime, and let $A_0,A_1, \ldots, A_k$ be
nonempty subsets of
the cyclic group $Z_p$.
If $|A_i| \neq |A_j|$ for all $0 \leq i<j \leq k$ and
$\sum_{i=0}^k |A_i| \leq p+{{k+2} \choose 2}-1$ then
$$
|\{a_0+a_1+ \ldots +a_k: a_i \in A_i,
a_i \neq a_j ~~\mbox{for all}~~i \neq j\}|
 \geq \sum_{i=0}^k |A_i|-{{k+2} \choose 2}+1.
$$
\end{prop}

The very special case of this proposition in which $k=1$,
$A_0=A$ and $A_1=A-\{a\}$ for an arbitrary element $a \in A$
implies that if $A \subset Z_p$ and $2|A|-1 \leq p+2$ then
the number of sums $a_1+a_2$ with $a_1,a_2 \in A$ and $a_1 \neq
a_2$ is at least $2|A|-3$. This supplies a short proof of
a result of Dias Da Silva and Hamidoune \cite{DH}, which settles a
conjecture of Erd\H os and Heilbronn
(cf., e.g., \cite{EGr}).

Snevily \cite{Sn} conjectured that for any two sets $A$ and $B$ of
equal cardinality in any
abelian group of odd order, there is a renumbering $a_i,b_i$ of
the  elements of $A$ and $B$ so that all sums
$a_i + b_i$ are pairwise
distinct.

For the cyclic group $Z_p$ of prime order, this follows
easily from Theorem \ref{t12} by considering the polynomial $f =
\prod_{i<j} (x_i-x_j) \prod_{i<j} (a_i+x_i-a_j-x_j)$
with $S_1 = \cdots = S_k = B$.

More generally, Dasgupta et al. \cite{DKSS} proved the conjecture for
any cyclic group of odd order, by applying the polynomial method for
polynomials over $Q[\omega]$, where $\omega$ is an appropriate root of
unity, and by considering $G$ as a subgroup of the multiplicative group
of this field. Further related results appear in \cite{Sun}.
%%%%%%%%%%

Additional applications  of Theorem \ref{t12} in additive number theory
can be found in \cite{Al3}.

\subsection{Graph coloring}

\vskip-5mm \hspace{5mm}

Theorem \ref{t12} has various applications in the study of
Graph Coloring, which is
the most popular area in Graph Theory. We sketch below
the basic approach, following
\cite{AT}. See also \cite{Mat}, \cite{Mat1} for a related
method.

A {\em vertex coloring} of a graph $G$ is an assignment of a color
to each vertex of $G$. The coloring is {\em proper} if
adjacent vertices  get distinct colors. The {\em chromatic
number}
$\chi(G)$ of $G$ is the minimum number of colors used in a proper
vertex coloring of $G$. An {\em edge coloring} of $G$ is,
similarly, an assignment of a color to each edge of $G$. It is
{\em proper} if adjacent edges receive distinct colors. The
minimum number of colors in a proper edge coloring of $G$ is
the {\em chromatic index} $\chi'(G)$ of $G$. This is
equal to the chromatic number of the line graph
of $G$.

A graph $G=(V,E)$ is $k$-{\em choosable} if for every assignment
of sets of integers $S(v) \subset Z$, each of size $k$, to the
vertices $v \in V$, there is a proper vertex coloring $c: V
\mapsto Z$ so that $c(v) \in S(v)$ for all $v \in V$. The {\em
choice number} of $G$, denoted by $ch(G)$, is the minimum integer
$k$ so that $G$ is $k$-choosable. Obviously, this number is at
least the chromatic number $\chi(G)$ of $G$. The choice number of
the line graph of $G$, denoted by $ch'(G)$, is usually called the
{\em list chromatic index} of $G$, and it is clearly at least the
chromatic index $\chi'(G)$ of $G$.

The study of choice numbers was
introduced, independently, by Vizing \cite{Vi} and by
Erd\H{o}s, Rubin and Taylor \cite{ERT}.
There are many graphs $G$ for which the choice number $ch(G)$ is
strictly larger than the chromatic number $\chi(G)$ (a complete
bipartite graph with $3$ vertices in each color class is one such
example).
In view of this, the following
conjecture, suggested independently by various researchers
including
Vizing, Albertson, Collins, Tucker and Gupta, which apparently
appeared first in print in
\cite{BH}, is somewhat surprising.
\begin{conj}{\rm( The list coloring conjecture)}
\label{c26}
For every graph $G$, $ch'(G)=\chi'(G)$.
\end{conj}

This conjecture asserts that for {\em line graphs} there is
no gap at all between the choice number and the chromatic number.
Many of the most interesting results in the area are proofs of
special
cases of this conjecture, which is still wide open.

The {\em graph polynomial}
$f_G=f_G(x_1,x_2, \ldots,x_n)$ of
a graph $G=(V,E)$ on
a set $V=\{1,\ldots,n\}$ of $n$ vertices is defined by
$f_G(x_1,x_2,\ldots ,x_n)=\Pi\big\{(x_i-x_j):i<j\ ,\
ij\in E\big\}$.
This polynomial has been studied by various researchers,
starting already with Petersen \cite{Pe} in 1891.

Note that if $S_1, \ldots ,S_n$ are sets of integers, then  there
is a proper coloring assigning to each vertex $i$ a color from its
list $S_i$, if and only if there are $s_i \in S_i$ such that
$f_G(s_1, \ldots ,s_n) \neq 0$. This condition is precisely the one
appearing in the conclusion of Theorem \ref{t12}, and it is
therefore natural to expect that this theorem can be useful in tackling
coloring problems. By applying it to line graphs of
planar cubic graphs, and by interpreting the appropriate
coefficient of the corresponding polynomial combinatorially, it can
be shown, using a known result of  Vigneron \cite{Vig} and
the Four Color Theorem, that the list chromatic index of
every $2$-connected cubic planar graph is $3$. This is
a strengthening of the Four Color
Theorem, which is well known to be equivalent to the fact that the
chromatic index of any such graph is $3$. An extension of this
result appears in \cite{EG}.

Additional results on graph coloring and choice numbers using the
above algebraic approach are described in the survey \cite{Al1}.
These include the fact that the choice number of every planar
bipartite graph is at most $3$, thus solving a conjecture raised
in \cite{ERT}, and the assertion, proved in \cite{FS}, that if $G$
is a graph on $3n$ vertices, whose set of edges is the disjoint
union of a Hamilton cycle and $n$ pairwise vertex-disjoint
triangles, then the choice number and the chromatic number of $G$
are both $3$.

\Section{The probabilistic method} \setzero

\vskip-5mm \hspace{5mm}

The discovery
that deterministic statements can be proved
by probabilistic reasoning, led already in the middle of the previous
century
to several striking results in Analysis,
Number Theory, Combinatorics and Information Theory.
It soon became clear that
the method, which is now called {\em the probabilistic method}, is
a very powerful tool for proving results in Discrete Mathematics.
The early results combined combinatorial arguments with fairly
elementary probabilistic techniques, whereas the development of the
method in recent years required the application of more
sophisticated tools from Probability Theory. In this section we
illustrate the method and describe several recent results.
More material can be found in the
recent books \cite{AS}, \cite{Bol}, \cite{JLR} and \cite{MRe}.

\subsection{Thresholds for random properties}

\vskip-5mm \hspace{5mm}

The systematic study of Random Graphs  was initiated by
Erd\H {o}s and R\'enyi whose first main paper on the subject is
\cite{ER}.
Formally, $G(n,p)$ denotes the probability space whose
points are graphs on a fixed set of $n$ labelled vertices, where
each pair of vertices forms an edge, randomly and independently,
with probability $p$. The term ``the random graph $G(n,p)$'' means,
in this context, a random point chosen in this probability space.
Each graph property $A$
(that is, a family of graphs closed under graph isomorphism)
is an event in this
probability space, and one may study its  probability $Pr[A]$,
that is, the probability that the random graph $G(n,p)$
lies in this family. In particular, we say that $A$ holds {\em
almost surely} if the probability that $G(n,p)$ satisfies $A$ tends
to $1$ as $n$ tends to infinity.
There are numerous papers dealing with random graphs, and
the two recent books \cite{Bol}, \cite{JLR} provide
excellent extensive accounts of
the known results in the subject.

One of the important discoveries of Erd\"os and R\'enyi
was the discovery of {\em threshold functions}. A function
$r(n)$ is called a threshold function for  a graph property $A$, if
when $p(n)/r(n)$ tends to $0$, then
$G(n,p(n))$ does not satisfy $A$ almost surely, whereas when
$p(n)/r(n)$  tends to infinity, then $G(n,p(n))$ satisfies $A$ almost
surely.  Thus, for
example, they identified the threshold function for the property
of being connected very precisely:
if $p(n)=\frac{\ln n}{n}+\frac{c}{n}$, then, as $n$ tends to
infinity, the probability that $G(n,p(n))$ is connected tends to
$e^{-e^{-c}}$.

A graph property is {\em monotone} if it is closed under the
addition  of edges.  Note that many interesting graph properties,
like hamiltonicity, non-planarity, connectivity or containing at
least $10$ vertex disjoint triangles
are monotone.

Bollob\'as and Thomason \cite{BT} proved that {\em any} monotone
graph property has a threshold function. Their proof applies to
any monotone family of subsets of a finite set, and holds
even without the assumption that the property $A$ is closed
under graph isomorphism.

Friedgut and Kalai \cite{FK} showed that the symmetry of graph
properties can be applied to obtain a sharper result. They proved
that for any
monotone graph property $A$, if $G(n,p)$ satisfies $A$
with probability at least $\epsilon$, then $G(n,q)$ satisfies $A$
with probability at least $1-\epsilon$, for
$q=p+O(\log (1/2 \epsilon) /\log n).$

The proof follows by
combining two results. The first is a
simple but fundamental lemma of Margulis \cite{Ma} and Russo
\cite{Ru}, which is useful in Percolation Theory. This lemma
can be used to
express the derivative with respect to $p$ of
the probability that $G(n,p)$ satisfies $A$  as a sum of
contributions associated with the single potential
edges. The second result is
a theorem of \cite{BKKKL}, which is proved using Harmonic Analysis,
that
asserts that at least one such contribution is always  large.
The symmetry implies that all contributions are the same and the result
follows. See also \cite{Ta} for
some related  results. These results hold for every transitive group
of symmetries.
In \cite{BK} it is shown that one can, in fact,
prove that the threshold for graph properties is even sharper, by
taking into account the precise group of symmetries induced
on the edges of the complete graph by permuting the vertices.
It turns  out that for every monotone graph property and for every
fixed $\epsilon>0$,
the width of the interval in which the probability
the property holds increases from $\epsilon$ to $1-\epsilon$ is at most
$c_{\delta}/(\log n)^{2-\delta}$ for all $\delta >0$. The power $2$
here is tight, as shown by the property of containing a clique of
size, say,  $\lfloor 2 \log_2 n \rfloor$.

It is natural to call the threshold for a monotone graph property
{\em sharp} if for every fixed positive $\epsilon$, the
width $w$ of the
interval in which the probability that the property holds
increases from $\epsilon$ to $1-\epsilon$ satisfies $w =o(p)$, where
$p$ is any point inside this interval. In \cite{Fr} Friedgut obtained
a beautiful characterization of all monotone
graph properties for which the threshold is sharp. Roughly speaking,
a property does not have a sharp threshold if and only if
it can be approximated
well in the relevant range of the probability $p$ by a property that
is determined by constant size witnesses. Thus, for example, the
property of containing $5$ vertex disjoint triangles does not
have a sharp threshold, whereas the property of having chromatic
number bigger than $10$ does. A similar result holds for  hypergraphs
as well. The proofs combine probabilistic and combinatorial
arguments with techniques from Harmonic analysis.

\subsection{Ramsey numbers}

\vskip-5mm \hspace{5mm}

Let $H_1, H_2, \ldots ,H_{k}$ be a
sequence of $k$ finite, undirected, simple
graphs. The (multicolored)
{\em Ramsey number} $r(H_1,H_2, \ldots ,H_{k})$ is the minimum integer
$r$ such that in every edge coloring
of the complete graph on $r$ vertices
by $k$ colors, there is a monochromatic copy of $H_i$ in color $i$ for
some $1 \leq i \leq k$. By a (special case of) a well known theorem of
Ramsey (c.f., e.g., \cite{GRS}), this number is finite for
every sequence of graphs $H_i$.

The determination or estimation of these numbers is usually a very
difficult problem.  When all graphs $H_i$ are complete graphs with
more than two vertices, the only values that are known precisely are
those of $r(K_3,K_m)$ for $m \leq 9$, $r(K_4,K_4)$, $r(K_4,K_5)$
and $r(K_3,K_3,K_3).$ Even the determination of the asymptotic
behaviour of Ramsey numbers  up to a constant factor is a hard
problem, and despite a lot of efforts by various researchers
(see, e.g., \cite{GRS}, \cite{CG} and their references), there
are only a few infinite families of graphs for which this behaviour
is known.

In one of the first applications of the probabilistic method in
Combinatorics, Erd\H{o}s \cite{Er1}
proved that if ${n \choose k}2^{1-{k \choose
2}}<1$ then $R(k,k)>n$, that is,
there exists a $2$-coloring of the edges of the complete graph
on $n$ vertices containing no monochromatic clique of
size $k$. The proof is
extremely simple;
the probability that a random two-edge coloring of $K_n$ contains
a monochromatic $K_k$ is at most
${n \choose k}2^{1-{k \choose 2}}<1$ , and hence there is
a coloring with the required  property.

A particularly interesting example of an infinite family
for which the asymtotic behaviour of the Ramsey number
is known,
is  the following result of Kim \cite{Ki} together
with that of Ajtai, Koml\'os and Szemer\'edi \cite{AKS0}.
\begin{theo}{\rm (\cite{Ki}, \cite{AKS0})}
\label{t81}
There are two absolute positive constants $c_1,c_2$ such that
$$ c_1 m^2 /\log m \leq
r(K_3,K_m) \leq c_2 m^2 / \log m $$
for all $m>1$.
\end{theo}

The upper bound, proved in \cite{AKS0}, is probabilistic, and applies
a certain random greedy algorithm. The lower bound is proved by
a ``semi-random''
construction and proceeds in stages. The detailed analysis
is subtle, and is based on certain large deviation inequalities.

Even less is known about
the asymptotic behaviour of multicolored Ramsey numbers,
that is, Ramsey numbers with at least $3$ colors.
The asymptotic
behaviour of $r(K_3,K_3,K_m)$, for example,
has been very poorly understood until recently, and
Erd\H{o}s and S\'os conjectured
in 1979 (c.f., e.g., \cite{CG}) that
$$
\lim_{m \mapsto \infty}  \frac{r(K_3,K_3,K_m)}{r(K_3,K_m)} =\infty.
$$
This has been proved recently, in a strong sense,
in \cite{AR}, where it is shown that
in fact $r(K_3,K_3,K_m)$ is equal, up to logarithmic factors,
to $m^3$.
A more complicated, related result  proved in \cite{AR}, that
supplies the asymptotic behaviour of infinitely many families of
Ramsey numbers up to a constant factor is the following.
\begin{theo}
\label{t82}
For every $t >1$ and $s \geq (t-1)!+1$ there are two
positive constants $c_1,c_2$ such that for every $m>1$
$$
c_1 \frac{m^t}{\log^t m} \leq r(K_{t,s},K_{t,s},K_{t,s},K_m)
\leq c_2 \frac{m^t}{\log^t m},
$$
where $K_{t,s}$ is the complete bipartite graph with $t$ vertices
in one color class and $s$ vertices in the other.
\end{theo}

The proof combines
spectral techniques, character sum estimates, and probabilistic
arguments.

\subsection{Tur\'an type results}

\vskip-5mm \hspace{5mm}

For a graph $H$ and an integer $n$, the
Tur\'an number $ex(n,H)$ is the maximum possible number of edges
in a simple graph on $n$ vertices that contains no copy of $H$.
The asymptotic behavior of these numbers for graphs of
chromatic number at least $3$ is well known, see, e.g.,
\cite{Bol0}. For bipartite graphs $H$, however, much less is known,
and there are relatively few nontrivial bipartite graphs $H$
for which the order of magnitude of $ex(n,H)$ is known.

A result of F\"uredi \cite{Fu} implies that for every fixed
bipartite graph $H$ in which the degrees of all vertices in one
color class are at most $r$, there is some $c=c(H)>0$ such that
$ex(n,H) \leq c n^{2-1/r}.$ As observed in \cite{AKS}, this result
can be derived from a simple and yet surprisingly powerful
probabilistic lemma, variants of which have been proved and
applied by various researchers starting with R\"odl and including
Kostochka, Gowers and Sudakov (see \cite{KR}, \cite{Go},
\cite{KS}). The lemma asserts, roughly, that every graph with
sufficiently many edges contains a large subset $A$ in which every
$a$ vertices have many common neighbors.  The proof uses a process
that may be called a {\em dependent random choice} for finding the
set $A$; $A$ is simply the set of all common neighbors of an
appropriately chosen random set $R$. Intuitively, it is clear that
if some $a$ vertices have only a few common neighbors, it is
unlikely all the members of $R$ will be chosen among these
neighbors. Hence, we do not expect $A$ to contain any such subset
of $a$ vertices. This simple idea can be extended. In particular,
it can be used to bound the Tur\'an numbers of degenerate
bipartite graphs.

A graph is {\em $r$-degenerate} if every subgraph of it contains
a vertex of degree at most $r$.  An old conjecture of Erd\H{o}s
asserts that for every fixed
$r$-degenerate bipartite graph $H$, $ex(n,H) \leq O(n^{2-1/r}),$
and the above technique suffices to show
that there is an absolute constant $c>0$, such that
for every such $H$, $ex(n,H) \leq n^{2-c/r}$.

Further questions and results about Tur\'an numbers  can be found in
\cite{AKS}, \cite{Bol0} and their references.

\Section{Algorithms and explicit constructions} \setzero

\vskip-5mm \hspace{5mm}

The rapid development of Theoretical Computer Science and its tight
connection to Discrete Mathematics motivated the study of the
algorithmic aspects of algebraic and probabilistic techniques. Can a
combinatorial structure, or a substructure of a given one,
whose existence is proved by algebraic or probabilistic
means, be constructed {\em explicitly} (that is, by an efficient
deterministic algorithm)? Can the algorithmic problems
corresponding to existence proofs be solved
by efficient procedures?  The study of these questions
often  requires tools from other
branches of mathematics.

As described in subsection 3.3, if $G$ is a graph on $3n$
vertices, whose set of edges is the disjoint union of a Hamilton
cycle and $n$ pairwise vertex-disjoint triangles, then the
chromatic number of $G$ is $3$. Can we solve the corresponding
algorithmic problem efficiently ? That is, is there a polynomial
time, deterministic or randomized algorithm, that given an input
graph as above, colors it properly with $3$ colors? Similarly, as
mentioned in subsection 3.3, the list chromatic index of any
planar cubic $2$-connected graph is $3$. Can we color properly the
edges of any given planar cubic $2$-connected graph using given
lists of three colors per edge, in polynomial time?

These problems, as well as the algorithmic versions of
additional applications of Theorem
\ref{t12}, are open. Of course, an algorithmic version
of the theorem itself would provide efficient procedures for solving all
these questions. The input for such  an algorithm
is a  polynomial in $n$ variables over a field
described, say, by a polynomial size arithmetic circuit.
Suppose that this polynomial
satisfies the assumptions of Theorem \ref{t12}, and that the
fact it satisfies it can be checked efficiently. The algorithm
should then find, efficiently, a point
$(s_1, s_2, \ldots ,s_n)$ satisfying the conclusion
of Theorem \ref{t12}.

Unfortunately, it seems unlikely that such a
general result can exist, as it would imply that there are no one-way
permutations. Indeed, let $F:\{0,1\}^n \mapsto \{0,1\}^n$ be a
$1-1$ function,
and suppose that for any $x=(x_1, \ldots ,x_n) \in
\{0,1\}^n$, the value of $F(x)$ can be computed efficiently.
Every Boolean function can be expressed as a multilinear
polynomial over $GF(2)$,
and hence, when we wish to find an  $x$ such that $F(x)=y=(y_1, \ldots
,y_n)$, we can write it as a system of multilinear
polynomials over $GF(2)$:
$F_i(x)=y_i$ for all $1 \leq i \leq n$. Equivalently,
this can be written
as $\prod_{i=1}^n (F_i(x)+y_i+1) \neq 0$.  This last equation
has a unique solution, implying that its left hand side, written as
a multilinear polynomial, is of full degree $n$ (since otherwise it
is easy to check that it attains the value $1$ an even number of times).
It follows that the assumptions of Theorem \ref{t12} with
$f=\prod_{i=1}^n (F_i(x)+y_i+1)$, $t_i=1$ and $S_i=GF(2)$ hold.
Thus, the existence of an efficient algorithm as above would
enable us to invert $F$ efficiently, implying that there cannot be any
one-way permutations. As this seems unlikely, it may be more
productive (and yet challenging)
to try and develop efficient procedures for solving
the particular algorithmic problems corresponding to the results
obtained by the theorem.

Probabilistic proofs also suggest the study of the corresponding
algorithmic problems. This is related to the study of randomized
algorithms, a topic which has been developed tremendously during
the last decade. See, e.g., \cite{MR} and its many references. In
particular, it is interesting to find explicit constructions of
combinatorial structures whose existence is proved by
probabilistic arguments. "Explicit" here means that there is a an
efficient algorithm that constructs the desired structure in time
polynomial in its size. Constructions of this type, besides being
interesting in their own, have applications in other areas. Thus,
for example, explicit constructions of error correcting codes that
are as good as the random ones are of interest in information
theory, and explicit constructions of certain Ramsey type
colorings may have applications in derandomization --- the process
of converting randomized algorithms into deterministic ones.

It turns out, however, that the problem of finding a good explicit
construction is often very difficult.
Even the simple proof of Erd\H os, described in
subsection 4.2,
that there are two-edge colorings of the complete graph on
$\lfloor 2^{m/2} \rfloor$
vertices containing no monochromatic clique
of size $m$, leads to an open
problem which seems very difficult. Can we construct, explicitly,
such a coloring of a complete graph on
$n \geq (1+\epsilon)^m$ vertices,
in time which is
polynomial in $n$, where $\epsilon>0$ is any positive absolute
constant ?

This problem is still open, despite a lot of efforts.
The best known explicit construction is due to Frankl and
Wilson \cite{FW}, who gave an explicit two-edge coloring of the
complete
graph on $m^{(1+o(1))\frac{\log m }{ 4 \log \log m}}$ vertices with
no monochromatic clique on $m$ vertices.

The construction of explicit two-edge colorings of large
complete graphs $K_n$
with no red $K_s$ and no blue $K_m$ for fixed $s$ and large $m$
also appears to be very difficult.
Using probabilistic arguments it can be shown that there are
such colorings for $n$ which is
$c \left(\frac{m}{\log m}\right)^{(s+1)/2}$ for some
absolute constant $c>0$. The best known explicit
construction, however, given in \cite{AP},
works only for
$m^{\delta \sqrt{\log s/\log \log s}}$, for some absolute
constant $\delta >0$. The description of the construction is
not complicated but the proof of its properties relies
on tools
from various mathematical areas. These include some ideas from
algebraic geometry obtained in \cite{KRS}, the well known bound of
Weil on character sums, spectral techniques and their connection to
the pseudo-random properties of graphs, the known bounds
of \cite{KST}
for the problem of Zarankiewicz  and the well known Erd\H{o}s-Rado
bound for the existence of $\Delta$-systems.

The above example is typical, and illustrates the fact that
tools from various mathematical disciplines often appear
in the design of explicit constructions of combinatorial
structures. Other examples that demonstrate this fact are the
construction of Algebraic Geometry codes, and the construction of
sparse pseudo-random graphs called expanders.

\Section{Some future challenges} \setzero

\vskip-5mm \hspace{5mm}

Several specific open problems in Discrete Mathematics are mentioned
throughout this article. These, and many additional ones, provide
interesting challenges for future research in the area.
We conclude with some brief comments on two more general
future challenges.

It seems safe to predict that in the future there will be
additional incorporation of methods from other mathematical areas
in Combinatorics. However, such methods often provide
non-constructive proof techniques, and the conversion of these to
algorithmic ones may well be one of the main future challenges of
the area. Another interesting recent development is the increased
appearance of Computer aided proofs in Combinatorics, starting
with the proof of the Four Color Theorem, and including automatic
methods for the discovery and proof of hypergeometric identities
--- see \cite{PWZ}. A successful incorporation of such proofs in
the area, without losing its special beauty and appeal, is another
challenge. These challenges, the fundamental nature of the area,
its tight connection to other disciplines, and the many
fascinating specific open problems studied in it, ensure that
Discrete Mathematics  will keep playing an essential role in the
general development of science in the future as well.

%\noindent
%{\bf Acknowledgment}

\label{lastpage}

\end{document}